# Hamacher and Heronian Aggregating operators under interval valued q-rung Orthopair Fuzzy and its Application

Xianping Li

**Abstract:** In this paper, combined with the degradation and monotonicity of parameters of Hamacher operation and the advantages of the correlation between input parameters of Heronian operators, IVq-ROFHHMWA and IVq-ROFHHMGA operators, which are the fusion of Hamacher and Heronian, are proposed, and their properties are studied. Finally, a group decision-making method based on IVq-ROFHHMWA operator is proposed. The implementation and comparative analysis results show that the operator and group decision-making methods proposed in this paper are effective and feasible.

## 1 Introduction

Based on Zadeh's fuzzy set theory [1], Atanassov's intuitionistic fuzzy set [2] takes into account the support, opposition and neutrality of decision makers, and has been widely used to solve multi-attribute decision-making problems in education, medical treatment and industry [3-5]. However, in the decision-making process, the sum of membership degree () and non-membership degree () must be less than or equal to 1, which increases the decision-making cost, so Yager proposed For uncertain fuzzy problems, researchers further put forward interval valued q-rung Orthopair Fuzzy sets [17-18], which have been well applied in group decision-making in recent years [9-13].

However, when solving group decision-making problems, aggregation operators are used to fuse the information of different experts and different attributes. Commonly used operators are: Haronian[19-21], Hamacher[15-16], Yager[22] and others [14]. Among them, Hamacher operation is characterized by Hamacher's t-conorms and t-norms, which are very popular in defining various fuzzy set operations and aggregation operators because of their advantages in parameter degeneration and monotonicity. Heronian mean operators (HM) are average polymerization technique, which is used to process accurate numerical values [23-24]. The ideal characteristic of HM is that it can capture the interrelation of input parameters, which makes it very useful for decision making [25]. Through investigation and analysis, so far, no scholars have studied and expanded Hamacher operation based on interval valued q-rung Orthopair Fuzzy, nor have they studied the aggregation operator of Hamacher and Haronian operator. Therefore, this paper proposes IVq-ROFHHMWA and IVq-ROFHHMGA operators in interval valued q-rung Orthopair Fuzzy environment, and studies their idempotent, monotonicity and permutation invariance. On the basis of IVq-ROFHHMWA operator, a new multi-attribute group decision-making method is proposed, which expands the degree of freedom of decision makers and has certain advantages in parameter degradation.

## 2 preliminary knowledge

### 2.1 Interval valued q-rung Orthopair Fuzzy sets

Definition 1[26-27]. Let x be the universe of discourse, and the definition of interval valued q-rung Orthopair Fuzzy set (IVq-ROFS)A on x is shown in formula (1).

$$A = \{(x, u_A(x), v_A(x)) | x \in X\} \quad (1)$$

In formula (1), the membership function is a mapping of interval values, satisfying: $u_A(x) = [u_A^-(x), u_A^+(x)] \subseteq [0,1]$, the non-membership function is a mapping of interval values, satisfying: $v_A(x) = [v_A^-(x), v_A^+(x)] \subseteq [0,1]$, simultaneously satisfying: $0 \le (u_A^+(x))^q + (v_A^+(x))^q \le 1, (q \ge 1)$, and its hesitation degree : $\pi_A(x) = [\pi_A^-(x), \pi_A^+(x)] = \left[\sqrt[q]{1-(u_A^+(x))^q-(v_A^+(x))^q},\ \sqrt[q]{1-(u_A^-(x))^q-(v_A^-(x))^q}\right]$

Definition 2[26-27]. Let $a = ([u^-, u^+], [v^-, v^+])$ $a_1 = ([u_{a_1}^-, u_{a_1}^+], [v_{a_1}^-, v_{a_1}^+])$ and $a_2 = ([u_{a_2}^-, u_{a_2}^+], [v_{a_2}^-, v_{a_2}^+])$ be three IVq-ROFNs and $q \ge 1$ where the operations of formulas (2), (3), (4) and (5) hold.

$$a_1 \oplus a_2 = <[\sqrt[q]{(u_{a_1}^-)^q + (u_{a_2}^-)^q - (u_{a_1}^-)^q (u_{a_2}^-)^q}, \sqrt[q]{(u_{a_1}^+)^q + (u_{a_2}^+)^q - (u_{a_1}^+)^q (u_{a_2}^+)^q}], [v_{a_1}^- v_{a_2}^-, v_{a_1}^+ v_{a_2}^+]> \quad (2)$$

$$a_1 \otimes a_2 = <[u_{a_1}^- u_{a_2}^-, u_{a_1}^+ u_{a_2}^+], [\sqrt[q]{(v_{a_1}^-)^q + (v_{a_2}^-)^q - (v_{a_1}^-)^q (v_{a_2}^-)^q}, \sqrt[q]{(v_{a_1}^+)^q + (v_{a_2}^+)^q - (v_{a_1}^+)^q (v_{a_2}^+)^q}]> \quad (3)$$

$$\lambda a = <[\sqrt[q]{1-(1-(u^-)^q)^\lambda}, \sqrt[q]{1-(1-(u^+)^q)^\lambda}], [(v^-)^\lambda, (v^+)^\lambda]> \quad (4)$$

$$a^\lambda = <[(u^-)^\lambda, (u^+)^\lambda], [\sqrt[q]{1-(1-(v^-)^q)^\lambda}, \sqrt[q]{1-(1-(v^+)^q)^\lambda}]> \quad (5)$$

Definition 3[28]. the score function of generalized orthogonal fuzzy numbers $a = < [u_a^-, u_a^+], [v_a^-, v_a^+] >$ with arbitrary interval values as shown in the following formula (6).

$$S(a) = \frac{1}{2}\left(u_a^- - v_a^+\left(1-u_a^+\right) + u_a^+ - v_a^-\left(1-u_a^-\right)\right) \quad (6)$$

Definition 4[26-27]. For any interval valued q-rung Orthopair Fuzzy number $a = < [u_a^-, u_a^+], [v_a^-, v_a^+] >$, its exact value function is shown in formula (7).

$$H(a) = \frac{1}{2}[(u_a^-)^q + (u_a^+)^q + (v_a^-)^q + (v_a^+)^q)], \quad (q \ge 1) \quad (7)$$

Definition 5[26-27]. For any two generalized orthogonal fuzzy numbers with interval values $a_1 = < [u_{a_1}^-, u_{a_1}^+], [v_{a_1}^-, v_{a_1}^+] >, a_2 = < [u_{a_2}^-, u_{a_2}^+], [v_{a_2}^-, v_{a_2}^+] >$, their magnitudes are compared as follows:

(1) If so $S(a_1) > S(a_2)$, then $a_1 > a_2$;

(2) If so $S(a_1) < S(a_2)$, then $a_1 < a_2$;

(3) If so $S(a_1) = S(a_2)$, further calculate their exact values and compare them. If $H(a_1) > H(a_2)$ ,$a_1 > a_2$; If $H(a_1) < H(a_2), a_1 < a_2$; If $H(a_1) = H(a_2), a_1 = a_2$。

### 2.2 Hamacher operation

Hamacher operation does not lose data in aggregation, and its clear domain is small. Although the calculation is somewhat complicated, the evaluation vector obtained is scientific and reasonable.

Definition 6[29] . For any two real numbers $a, b \in [0, 1]$, Hamacher sum $\oplus_H$ is a T-norm operation, and Hamacher product $\otimes_H$ is a T-norm operation, such as formulas (8) and (9).

$$a \oplus_H b = \frac{a+b-ab-(1-\varphi)ab}{1-(1-\varphi)(ab)}, \varphi > 0 \tag{8}$$

$$a \otimes_H b = \frac{ab}{\varphi+(1-\varphi)(a+b-ab)}, \varphi > 0 \tag{9}$$

## 2.3 Heronian mean operators

Heronian mean operators (HM) is an excellent operator considering the interaction of aggregation variables, and the changes of x and y under parameters can reflect the sensitivity of aggregation.

Definition 7[30]. Let $x, y \geq 0$, and not be 0 at the same time, $a_i (i = 1,2, \dots n)$ be a set of real numbers. If HM: $HM^n \rightarrow HM$ is a mapping, Heronian mean operators is shown in formula (10).

$$HM^{x,y}\left(a_{1,}a_{2,}a_3, \dots a_n\right) = \left(\frac{2}{n(n+1)}\left(\sum_{i=1,j=i}^{n}\left(a_i^x a_j^y\right)\right)\right)^{\frac{1}{x+y}} \tag{10}$$

# 3 IVq-ROFHHM operator

## 3.1 IVq-ROFHHM algorithm

According to t-conorms and t-norms of Hamacher operation in formulas (8) and (9), this paper defines the arithmetic rules of interval valued q-rung Orthopair Fuzzy Hamacher addition, multiplication, number multiplication and power.

Definition 8. Let $a = [u^-, u^+], [v^-, v^+], a_1 = [u_1^-, u_1^+], [v_1^-, v_1^+]$, $a_2 = [u_2^-, u_2^+], [v_2^-, v_2^+]$ be three q-RIVOFNs,the addition, multiplication, number multiplication and power operation of interval valued q-rung Orthopair Fuzzy Hamacher operation refer to formulas (11)-(14) directly.

$$a_{1\oplus_H}a_2 = \left(\begin{matrix}\left[\begin{matrix}\sqrt[q]{\frac{(u_1^-)^q+(u_2^-)^q-(u_1^-)^q(u_2^-)^q-(1-\varphi)(u_1^-)^q(u_2^-)^q}{1-(1-\varphi)(u_1^-)^q(u_2^-)^q}}, \\ \sqrt[q]{\frac{(u_1^+)^q+(u_2^+)^q-(u_1^+)^q(u_2^+)^q-(1-\varphi)(u_1^+)^q(u_2^+)^q}{1-(1-\varphi)(u_1^+)^q(u_2^+)^q}}\end{matrix}\right], \\ \left[\begin{matrix}\frac{v_1^- v_2^-}{\sqrt[q]{\varphi+(1-\varphi)((v_1^-)^q+(v_2^-)^q-(v_1^-)^q(v_2^-)^q)}}, \\ \frac{v_1^+ v_2^+}{\sqrt[q]{\varphi+(1-\varphi)\left((v_1^+)^q+(v_2^+)^q-(v_1^+)^q(v_2^+)^q\right)}}\end{matrix}\right]\end{matrix}\right), \varphi > 0 \tag{11}$$

$$a_{1\otimes_H}a_2 = \left(\begin{matrix}\left[\begin{matrix}\frac{u_1^- u_2^-}{\sqrt[q]{\varphi+(1-\varphi)((u_1^-)^q+(u_2^-)^q-(u_1^-)^q(u_2^-)^q)}}, \\ \frac{u_1^+ u_2^+}{\sqrt[q]{\varphi+(1-\varphi)\left((u_1^+)^q+(u_2^+)^q-(u_1^+)^q(u_2^+)^q\right)}}\end{matrix}\right], \\ \left[\begin{matrix}\sqrt[q]{\frac{(v_1^-)^q+(v_2^-)^q-(v_1^-)^q(v_2^-)^q-(1-\varphi)(v_1^-)^q(v_2^-)^q}{1-(1-\varphi)(v_1^-)^q(v_2^-)^q}}, \\ \sqrt[q]{\frac{(v_1^+)^q+(v_2^+)^q-(v_1^+)^q(v_2^+)^q-(1-\varphi)(v_1^+)^q(v_2^+)^q}{1-(1-\varphi)(v_1^+)^q(v_2^+)^q}}\end{matrix}\right]\end{matrix}\right), \varphi > 0 \tag{12}$$

$$\vartheta a=\left(\begin{array}{c}\left[\sqrt[q]{\frac{(1+(\varphi-1)(u^-)^q)^\vartheta-(1-(u^-)^q)^\vartheta}{(1+(\varphi-1)(u^-)^q)^\vartheta+(\varphi-1)(1-(u^-)^q)^\vartheta}},\sqrt[q]{\frac{(1+(\varphi-1)(u^+)^q)^\vartheta-(1-(u^+)^q)^\vartheta}{(1+(\varphi-1)(u^+)^q)^\vartheta+(\varphi-1)(1-(u^+)^q)^\vartheta}}\right],\\ \left[\frac{\sqrt[q]{\varphi}(v^-)^\vartheta}{\sqrt[q]{\left(1+(\varphi-1)(1-(v^-)^q)\right)^\vartheta+(\varphi-1)((v^-)^q)^\vartheta}},\frac{\sqrt[q]{\varphi}(v^+)^\vartheta}{\sqrt[q]{\left(1+(\varphi-1)(1-(v^+)^q)\right)^\vartheta+(\varphi-1)((v^+)^q)^\vartheta}}\right]\end{array}\right),\varphi>0 \quad (13)$$

$$a^\vartheta=\left(\begin{array}{c}\left[\frac{\sqrt[q]{\varphi}(u^-)^\vartheta}{\sqrt[q]{\left(1+(\varphi-1)(1-(u^-)^q)\right)^\vartheta+(\varphi-1)((u^-)^q)^\vartheta}},\frac{\sqrt[q]{\varphi}(u^+)^\vartheta}{\sqrt[q]{\left(1+(\varphi-1)(1-(u^+)^q)\right)^\vartheta+(\varphi-1)((u^+)^q)^\vartheta}}\right],\\ \left[\sqrt[q]{\frac{(1+(\varphi-1)(v^-)^q)^\vartheta-(1-(v^-)^q)^\vartheta}{(1+(\varphi-1)(v^-)^q)^\vartheta+(\varphi-1)(1-(v^-)^q)^\vartheta}},\sqrt[q]{\frac{(1+(\varphi-1)(v^+)^q)^\vartheta-(1-(v^+)^q)^\vartheta}{(1+(\varphi-1)(v^+)^q)^\vartheta+(\varphi-1)(1-(v^+)^q)^\vartheta}}\right]\end{array}\right),\varphi>0 \quad (14)$$

It can be verified that the addition and multiplication of interval valued q-rung Orthopair Fuzzy Hamacher operation given in this paper satisfy the operation rules of t-conorms and t-norms, and the multiplication and power of numbers also satisfy the operation rules of t-conorms and T-conorms.

## 3.2 IVq-ROFHHM operator

Definition 9. Let $x, y \geq 0$, not be 0 at the same time , $a_i = \langle[u_i^-, u_i^+], [v_i^-, v_i^+]\rangle (i = 1,2, \ldots n)$ is a group of IVq-ROFNs. There is a mapping: $IVq-ROFHHM^n \rightarrow IVq-ROFHHM$, then the IVq-ROFHHM operator is shown in Formula (15).

$$IVq-ROFHHM^{x,y}(a_1,a_2,\ldots a_n)=\frac{1}{x+y}\left(\underset{i=1,j=1}{\overset{n}{\otimes_H}}\left(xa_i\oplus_H ya_j\right)^{\frac{2}{n(n+1)}}\right) \quad (15)$$

Theorem 1: Let $x, y \geq 0$, not be 0 at the same time and $a_i = \langle[u_i^-, u_i^+], [v_i^-, v_i^+]\rangle$ $(i = 1,2, \ldots n)$ is a group of IVq-ROFNs, and the operation formula of IVq-ROFHHM operator is shown in (16).

$$IVq-ROFHHM^{x,y}(a_1,a_2,\ldots a_n)=\frac{1}{x+y}\left(\underset{i=1,j=1}{\overset{n}{\otimes_H}}\left(xa_i\oplus_H ya_j\right)^{\frac{2}{n(n+1)}}\right)$$

$$=\left(\begin{array}{l}\left[\sqrt[q]{\frac{a^- -b^-}{a^-+(\varphi-1)b^-}},\sqrt[q]{\frac{a^+ -b^+}{a^++(\varphi-1)b^+}}\right],\\ \left[\sqrt[q]{\frac{\varphi\left(c^--d^-\right)^{\frac{1}{x+y}}}{\left(c^-+(\varphi^2-1)d^-\right)^{\frac{1}{x+y}}+(\varphi-1)\left(c^--d^-\right)^{\frac{1}{x+y}}}},\sqrt[q]{\frac{\varphi\left(c^+-d^+\right)^{\frac{1}{x+y}}}{\left(c^++(\varphi^2-1)d^+\right)^{\frac{1}{x+y}}+(\varphi-1)\left(c^+-d^+\right)^{\frac{1}{x+y}}}}\right]\end{array}\right) \quad (16)$$

In formula (16), there are:

$$a^-=\left(\prod_{i=1,j=1}^{n}\left(V_{ij}^-+(\varphi^2-1)W_{ij}^-\right)^{\frac{2}{n(n+1)}}+(\varphi^2-1)\prod_{i=1,j=1}^{n}\left(V_{ij}^--W_{ij}^-\right)^{\frac{2}{n(n+1)}}\right)^{\frac{1}{x+y}}$$

$$a^{+}=\left(\prod_{i=1,j=1}^{n}\left(V_{ij}^{+}+\left(\varphi^{2}-1\right)W_{ij}^{+}\right)^{\frac{2}{n(n+1)}}+\left(\varphi^{2}-1\right)\prod_{i=1,j=1}^{n}\left(V_{ij}^{+}-W_{ij}^{+}\right)^{\frac{2}{n(n+1)}}\right)^{\frac{1}{x+y}}$$

$$b^{-}=\left(\prod_{i=1,j=1}^{n}\left(V_{ij}^{-}+\left(\varphi^{2}-1\right)W_{ij}^{-}\right)^{\frac{2}{n(n+1)}}-\prod_{i=1,j=1}^{n}\left(V_{ij}^{-}-W_{ij}^{-}\right)^{\frac{2}{n(n+1)}}\right)^{\frac{1}{x+y}}$$

$$b^{+}=\left(\prod_{i=1,j=1}^{n}\left(V_{ij}^{+}+\left(\varphi^{2}-1\right)W_{ij}^{+}\right)^{\frac{2}{n(n+1)}}-\prod_{i=1,j=1}^{n}\left(V_{ij}^{+}-W_{ij}^{+}\right)^{\frac{2}{n(n+1)}}\right)^{\frac{1}{x+y}}$$

$$c^{-}=\left(\prod_{i=1,j=1}^{n}\left(N_{ij}^{-}+\left(\varphi^{2}-1\right)M_{ij}^{-}\right)^{\frac{2}{n(n+1)}}-\prod_{i=1,j=1}^{n}\left(N_{ij}^{-}-M_{ij}^{-}\right)^{\frac{2}{n(n+1)}}\right)^{\frac{1}{x+y}}$$

$$c^{+}=\left(\prod_{i=1,j=1}^{n}\left(N_{ij}^{+}+\left(\varphi^{2}-1\right)M_{ij}^{+}\right)^{\frac{2}{n(n+1)}}-\prod_{i=1,j=1}^{n}\left(N_{ij}^{+}-M_{ij}^{+}\right)^{\frac{2}{n(n+1)}}\right)^{\frac{1}{x+y}}$$

$$d^{-}=\left(\prod_{i=1,j=1}^{n}\left(N_{ij}^{-}+\left(\varphi^{2}-1\right)M_{ij}^{-}\right)^{\frac{2}{n(n+1)}}-\left(\varphi^{2}-1\right)\prod_{i=1,j=1}^{n}\left(N_{ij}^{-}-M_{ij}^{-}\right)^{\frac{2}{n(n+1)}}\right)^{\frac{1}{x+y}}$$

$$+(\varphi-1)\left(\prod_{i=1,j=1}^{n}\left(N_{ij}^{-}+\left(\varphi^{2}-1\right)M_{ij}^{-}\right)^{\frac{2}{n(n+1)}}-\prod_{i=1,j=1}^{n}\left(N_{ij}^{-}-M_{ij}^{-}\right)^{\frac{2}{n(n+1)}}\right)^{\frac{1}{x+y}}$$

$$d^{+}=\left(\prod_{i=1,j=1}^{n}\left(N_{ij}^{+}+\left(\varphi^{2}-1\right)M_{ij}^{+}\right)^{\frac{2}{n(n+1)}}-\left(\varphi^{2}-1\right)\prod_{i=1,j=1}^{n}\left(N_{ij}^{+}-M_{ij}^{+}\right)^{\frac{2}{n(n+1)}}\right)^{\frac{1}{x+y}}$$

$$+(\varphi-1)\left(\prod_{i=1,j=1}^{n}\left(N_{ij}^{+}+\left(\varphi^{2}-1\right)M_{ij}^{+}\right)^{\frac{2}{n(n+1)}}-\prod_{i=1,j=1}^{n}\left(N_{ij}^{+}-M_{ij}^{+}\right)^{\frac{2}{n(n+1)}}\right)^{\frac{1}{x+y}}$$

Where in the above formula satisfies the following formula:

$$V_{ij}^{-}=\left(1+(\varphi-1)\left(1-\left(\mu_{i}^{-}\right)^{q}\right)+\left(\varphi^{2}-1\right)\left(\mu_{i}^{-}\right)^{q}\right)^{x}\left(1+(\varphi-1)\left(1-\left(\mu_{j}^{-}\right)^{q}\right)+\left(\varphi^{2}-1\right)\left(\mu_{j}^{-}\right)^{q}\right)^{y}$$

$$V_{ij}^{+}=\left(1+(\varphi-1)\left(1-\left(\mu_{i}^{+}\right)^{q}\right)+\left(\varphi^{2}-1\right)\left(\mu_{i}^{+}\right)^{q}\right)^{x}\left(1+(\varphi-1)\left(1-\left(\mu_{j}^{+}\right)^{q}\right)+\left(\varphi^{2}-1\right)\left(\mu_{j}^{+}\right)^{q}\right)^{y}$$

$$W_{ij}^{-}=\left(1+(\varphi-1)\left(1-\left(\mu_{i}^{-}\right)^{q}\right)-\left(\mu_{i}^{-}\right)^{q}\right)^{x}\left(1+(\varphi-1)\left(1-\left(\mu_{j}^{-}\right)^{q}\right)-\left(\mu_{j}^{-}\right)^{q}\right)^{y}$$

$$W_{ij}^{+}=\left(1+(\varphi-1)\left(1-\left(\mu_{i}^{+}\right)^{q}\right)-\left(\mu_{i}^{+}\right)^{q}\right)^{x}\left(1+(\varphi-1)\left(1-\left(\mu_{j}^{+}\right)^{q}\right)-\left(\mu_{j}^{+}\right)^{q}\right)^{y}$$

$$N_{ij}^{-}=\left(1+(\varphi-1)\left(1-\left(v_{i}^{-}\right)^{q}\right)-\left(\varphi^{2}-1\right)\left(1-\left(v_{i}^{-}\right)^{q}\right)\right)^{x}\left(1+(\varphi-1)\left(v_{j}^{-}\right)^{q}-\left(\varphi^{2}-1\right)\left(1-\left(v_{j}^{-}\right)^{q}\right)\right)^{y}$$

$$N_{ij}^{+}=\left(1+(\varphi-1)\left(1-\left(v_{i}^{+}\right)^{q}\right)-\left(\varphi^{2}-1\right)\left(1-\left(v_{i}^{+}\right)^{q}\right)\right)^{x}\left(1+(\varphi-1)\left(v_{j}^{+}\right)^{q}-\left(\varphi^{2}-1\right)\left(1-\left(v_{j}^{+}\right)^{q}\right)\right)^{y}$$

$$M_{ij}^{-}=\left(\varphi\left(v_i^{-}\right)^q\right)^x\left(\varphi\left(v_j^{-}\right)^q\right)^y$$

$$M_{ij}^{-}=\left(\varphi\left(v_i^{+}\right)^q\right)^x\left(\varphi\left(v_j^{+}\right)^q\right)^y$$

Combining with the interval valued q-rung Orthopair Fuzzy HHM operator algorithms (11)-(14), the formula (16) is proved by induction, namely

When n=2 $IVq-ROFHHM^{x,y}\left(a_1,a_2\right)=\frac{1}{x+y}\left(\left(xa_1\oplus_H ya_2\right)\otimes_H\left(ya_1\oplus_H xa_2\right)\right)^{\frac{1}{3}}$

Let n=k-1 $IVq-ROFHHM^{x,y}\left(a_1,a_2,...a_{k-1}\right)=\frac{1}{x+y}\left(\underset{i=1,j=1}{\overset{k-1}{\otimes_H}}\left(xa_i\oplus_H ya_j\right)^{\frac{2}{k(k-1)}}\right)$

When n=k

$$IVq-ROFHHM^{x,y}\left(a_1,a_2,...a_k\right)$$

$$=\frac{1}{x+y}\left(\left(\left(x\left(\underset{i=1,j=1}{\overset{k-1}{\otimes_H}}\left(xa_i\oplus_H ya_j\right)^{\frac{2}{k(k-1)}}\right)\oplus_H ya_k\right)\otimes_H\left(xa_k\oplus_H y\left(\underset{i=1,j=1}{\overset{k-1}{\otimes_H}}\left(xa_i\oplus_H ya_j\right)^{\frac{2}{k(k-1)}}\right)\right)\right)\right)^{\frac{1}{3}}$$

$$=\left(\begin{array}{l}\left[\sqrt[q]{\frac{a^{-}-b^{-}}{a^{-}+(\varphi-1)b^{-}}},\sqrt[q]{\frac{a^{+}-b^{+}}{a^{+}+(\varphi-1)b^{+}}}\right],\\ \left[\sqrt[q]{\frac{\varphi\left(c^{-}-d^{-}\right)^{\frac{1}{x+y}}}{\left(c^{-}+\left(\varphi^2-1\right)d^{-}\right)^{\frac{1}{x+y}}+(\varphi-1)\left(c^{-}-d^{-}\right)^{\frac{1}{x+y}}}},\sqrt[q]{\frac{\varphi\left(c^{+}-d^{+}\right)^{\frac{1}{x+y}}}{\left(c^{+}+\left(\varphi^2-1\right)d^{+}\right)^{\frac{1}{x+y}}+(\varphi-1)\left(c^{+}-d^{+}\right)^{\frac{1}{x+y}}}}\right]\end{array}\right)\otimes_H a_k$$

and

$$a^{-}=\left(\prod_{i=1,j=1}^{k-1}\left(V_{ij}^{-}+\left(\varphi^2-1\right)W_{ij}^{-}\right)^{\frac{2}{k(k-1)}}+\left(\varphi^2-1\right)\prod_{i=1,j=1}^{k-1}\left(V_{ij}^{-}-W_{ij}^{-}\right)^{\frac{2}{k(k-1)}}\right)^{\frac{1}{x+y}}$$

$$a^{+}=\left(\prod_{i=1,j=1}^{k-1}\left(V_{ij}^{+}+\left(\varphi^2-1\right)W_{ij}^{+}\right)^{\frac{2}{k(k-1)}}+\left(\varphi^2-1\right)\prod_{i=1,j=1}^{k-1}\left(V_{ij}^{+}-W_{ij}^{+}\right)^{\frac{2}{k(k-1)}}\right)^{\frac{1}{x+y}}$$

$$b^{-}=\left(\prod_{i=1,j=1}^{k-1}\left(V_{ij}^{-}+\left(\varphi^2-1\right)W_{ij}^{-}\right)^{\frac{2}{k(k-1)}}-\prod_{i=1,j=1}^{k-1}\left(V_{ij}^{-}-W_{ij}^{-}\right)^{\frac{2}{k(k-1)}}\right)^{\frac{1}{x+y}}$$

$$b^{+}=\left(\prod_{i=1,j=1}^{k-1}\left(V_{ij}^{+}+\left(\varphi^2-1\right)W_{ij}^{+}\right)^{\frac{2}{k(k-1)}}-\prod_{i=1,j=1}^{k-1}\left(V_{ij}^{+}-W_{ij}^{+}\right)^{\frac{2}{k(k-1)}}\right)^{\frac{1}{x+y}}$$

$$c^{-}=\left(\prod_{i=1,j=1}^{k-1}\left(N_{ij}^{-}+\left(\varphi^2-1\right)M_{ij}^{-}\right)^{\frac{2}{k(k-1)}}-\prod_{i=1,j=1}^{k-1}\left(N_{ij}^{-}-M_{ij}^{-}\right)^{\frac{2}{k(k-1)}}\right)^{\frac{1}{x+y}}$$

$$c^+ = \left( \prod_{i=1,\, j=1}^{k-1} \left( N_{ij}^+ + \left(\varphi^2 - 1\right) M_{ij}^+ \right)^{\frac{2}{k(k-1)}} - \prod_{i=1,\, j=1}^{k-1} \left( N_{ij}^+ - M_{ij}^+ \right)^{\frac{2}{k(k-1)}} \right)^{\frac{1}{x+y}}$$

$$d^- = \left( \prod_{i=1,\, j=1}^{k-1} \left( N_{ij}^- + \left(\varphi^2 - 1\right) M_{ij}^- \right)^{\frac{2}{k(k-1)}} - \left(\varphi^2 - 1\right) \prod_{i=1,\, j=1}^{k-1} \left( N_{ij}^- - M_{ij}^- \right)^{\frac{2}{k(k-1)}} \right)^{\frac{1}{x+y}}$$

$$+ \left(\varphi - 1\right) \left( \prod_{i=1,\, j=1}^{k-1} \left( N_{ij}^- + \left(\varphi^2 - 1\right) M_{ij}^- \right)^{\frac{2}{k(k-1)}} - \prod_{i=1,\, j=1}^{k-1} \left( N_{ij}^- - M_{ij}^- \right)^{\frac{2}{k(k-1)}} \right)^{\frac{1}{x+y}}$$

$$d^+ = \left( \prod_{i=1,\, j=1}^{k-1} \left( N_{ij}^+ + \left(\varphi^2 - 1\right) M_{ij}^+ \right)^{\frac{2}{k(k-1)}} - \left(\varphi^2 - 1\right) \prod_{i=1,\, j=1}^{k-1} \left( N_{ij}^+ - M_{ij}^+ \right)^{\frac{2}{k(k-1)}} \right)^{\frac{1}{x+y}}$$

$$+ \left(\varphi - 1\right) \left( \prod_{i=1,\, j=1}^{k-1} \left( N_{ij}^+ + \left(\varphi^2 - 1\right) M_{ij}^+ \right)^{\frac{2}{k(k-1)}} - \prod_{i=1,\, j=1}^{k-1} \left( N_{ij}^+ - M_{ij}^+ \right)^{\frac{2}{k(k-1)}} \right)^{\frac{1}{x+y}}$$

$$= \left( \begin{array}{l} \left[ \sqrt[q]{\dfrac{a^- - b^-}{a^- + (\varphi - 1) b^-}}, \sqrt[q]{\dfrac{a^+ - b^+}{a^+ + (\varphi - 1) b^+}} \right], \\ \left[ \sqrt[q]{\dfrac{\varphi \left(c^- - d^-\right)^{\frac{1}{x+y}}}{\left(c^- + \left(\varphi^2 - 1\right) d^-\right)^{\frac{1}{x+y}} + (\varphi - 1)\left(c^- - d^-\right)^{\frac{1}{x+y}}}}, \sqrt[q]{\dfrac{\varphi \left(c^+ - d^+\right)^{\frac{1}{x+y}}}{\left(c^+ + \left(\varphi^2 - 1\right) d^+\right)^{\frac{1}{x+y}} + (\varphi - 1)\left(c^+ - d^+\right)^{\frac{1}{x+y}}}} \right] \end{array} \right), \text{ and}$$

$$a^- = \left( \prod_{i=1,\, j=1}^{k} \left( V_{ij}^- + \left(\varphi^2 - 1\right) W_{ij}^- \right)^{\frac{2}{k(k+1)}} + \left(\varphi^2 - 1\right) \prod_{i=1,\, j=1}^{k} \left( V_{ij}^- - W_{ij}^- \right)^{\frac{2}{k(k+1)}} \right)^{\frac{1}{x+y}}$$

$$a^+ = \left( \prod_{i=1,\, j=1}^{k} \left( V_{ij}^+ + \left(\varphi^2 - 1\right) W_{ij}^+ \right)^{\frac{2}{k(k+1)}} + \left(\varphi^2 - 1\right) \prod_{i=1,\, j=1}^{k} \left( V_{ij}^+ - W_{ij}^+ \right)^{\frac{2}{k(k+1)}} \right)^{\frac{1}{x+y}}$$

$$b^- = \left( \prod_{i=1,\, j=1}^{k} \left( V_{ij}^- + \left(\varphi^2 - 1\right) W_{ij}^- \right)^{\frac{2}{k(k+1)}} - \prod_{i=1,\, j=1}^{k} \left( V_{ij}^- - W_{ij}^- \right)^{\frac{2}{k(k+1)}} \right)^{\frac{1}{x+y}}$$

$$b^+ = \left( \prod_{i=1,\, j=1}^{k} \left( V_{ij}^+ + \left(\varphi^2 - 1\right) W_{ij}^+ \right)^{\frac{2}{k(k+1)}} - \prod_{i=1,\, j=1}^{k} \left( V_{ij}^+ - W_{ij}^+ \right)^{\frac{2}{k(k+1)}} \right)^{\frac{1}{x+y}}$$

$$c^- = \left( \prod_{i=1,\, j=1}^{k} \left( N_{ij}^- + \left(\varphi^2 - 1\right) M_{ij}^- \right)^{\frac{2}{k(k+1)}} - \prod_{i=1,\, j=1}^{k} \left( N_{ij}^- - M_{ij}^- \right)^{\frac{2}{k(k+1)}} \right)^{\frac{1}{x+y}}$$

$$c^+ = \left( \prod_{i=1,\, j=1}^{k} \left( N_{ij}^+ + \left(\varphi^2 - 1\right) M_{ij}^+ \right)^{\frac{2}{k(k+1)}} - \prod_{i=1,\, j=1}^{k} \left( N_{ij}^+ - M_{ij}^+ \right)^{\frac{2}{k(k+1)}} \right)^{\frac{1}{x+y}}$$

$$d^- = \left( \prod_{i=1,\, j=1}^{k} \left( N_{ij}^- + \left(\varphi^2 - 1\right) M_{ij}^- \right)^{\frac{2}{k(k+1)}} - \left(\varphi^2 - 1\right) \prod_{i=1,\, j=1}^{k} \left( N_{ij}^- - M_{ij}^- \right)^{\frac{2}{k(k+1)}} \right)^{\frac{1}{x+y}}$$

$$+(\varphi-1)\left(\prod_{i=1,j=1}^{k}\left(N_{ij}^{-}+(\varphi^{2}-1)M_{ij}^{-}\right)^{\frac{2}{k(k+1)}}-\prod_{i=1,j=1}^{k}\left(N_{ij}^{-}-M_{ij}^{-}\right)^{\frac{2}{k(k+1)}}\right)^{\frac{1}{x+y}}$$

$$d^{+}=\left(\prod_{i=1,j=1}^{k}\left(N_{ij}^{+}+(\varphi^{2}-1)M_{ij}^{+}\right)^{\frac{2}{k(k+1)}}-(\varphi^{2}-1)\prod_{i=1,j=1}^{k}\left(N_{ij}^{+}-M_{ij}^{+}\right)^{\frac{2}{k(k+1)}}\right)^{\frac{1}{x+y}}$$

$$+(\varphi-1)\left(\prod_{i=1,j=1}^{k}\left(N_{ij}^{+}+(\varphi^{2}-1)M_{ij}^{+}\right)^{\frac{2}{k(k+1)}}-\prod_{i=1,j=1}^{k}\left(N_{ij}^{+}-M_{ij}^{+}\right)^{\frac{2}{k(k+1)}}\right)^{\frac{1}{x+y}}$$

Prove the theorem.

IVq-ROFHHM operator satisfies permutation, boundary, monotonicity. Below, we give the properties of IVq-ROFHHM operator.

Theorem 2 (Permutation): If there is any permutation, then IVq-ROFHHM operator satisfies permutation, as shown in formula (17).$(\beta_1,\beta_2,\dots\beta_n)(a_1,a_2,\dots a_n)$

$$IVq-ROFHHM^{x,y}(\beta_1,\beta_2,\dots\beta_n)=IVq-ROFHHM^{x,y}(a_1,a_2,\dots a_n)\quad(17)$$

Prove:

Let $(a_1^{'},a_2^{'},\dots,a_n^{'})$ be any permutation of $(a_1,a_2,\dots,a_n)$, and for any $a_i^{'}$, it corresponds to unique $a_j$. So:

$$IVq-ROFHHM^{x,y}(a_1,a_2,\dots a_n)=\frac{1}{x+y}\left(\underset{i=1,j=1}{\overset{n}{\otimes_H}}\left(xa_i\oplus_H ya_j\right)^{\frac{2}{n(n+1)}}\right)$$

$$=IVq-ROFHHM^{x,y}(a_1^{'},a_2^{'},\dots,a_n^{'})$$

Theorem 3 (Monotonicity): Let $a_i=\langle[u_i^-,u_i^+],[v_i^-,v_i^+]\rangle$, $a_i'=\langle[u_i^{-\prime},u_i^{+\prime}],[v_i^{-\prime},v_i^{+\prime}]\rangle(i=1,2,\dots,\text{n})$ are two groups of IVq-ROFNs and satisfy: $u_i^-\le u_i^{-\prime}$, $u_i^+\le u_i^{+\prime}$, $v_i^-\ge v_i^{-\prime}$, $v_i^+\ge v_i^{+\prime}$, then the IVq-ROFHHM operator satisfies monotonicity, as shown in formula (18).

$$IVq-ROFHHM^{x,y}(a_1,a_2,\dots a_n)\le IVq-ROFHHM^{x,y}(a_1',a_2'\dots a_n')\quad(18)$$

Prove:

Substituting the formula (16) into the scoring function formula (6) to calculate: $S(IVq-ROFHHM(a_1,a_2,\dots,a_n))$, $S(IVq-ROFHHM(a_1',a_2'\dots a_n'))$

Because $u_{a_i}^-\le u_{a_i^{'}}^-,u_{a_i}^+\le u_{a_i^{'}}^+,\ v_{a_i}^-\ge v_{a_i^{'}}^-,v_{a_i}^+\ge v_{a_i^{'}}^+$

Conclude $S(IVq-ROFHHM^{x,y}(a_1,a_2,\dots,a_n))\le S\left(IVq-ROFHHM^{x,y}(a_1^{'},a_2^{'}\dots a_n^{'})\right)$

According to the size comparison definition, you can get:

$$IVq-ROFHHM^{x,y}(a_1,a_2,\dots a_n)\le IVq-ROFHHM^{x,y}(a_1^{'},a_2^{'}\dots a_n^{'})$$

It is worth noting that the monotonicity given here needs to satisfy the size correspondence relationship of the corresponding interval-valued generalized orthogonal fuzzy. When the data does not satisfy this relationship, the monotonicity may not necessarily hold.

Theorem 4 (Boundary): Let $(a_1,a_2,\dots a_n)$ be a set of IVq-ROFNs, $a_{\min}=\min_i\{a_i\}$, $a_{\max}=\max_i\{a_i\}$. Then the IVq-ROFHHM operator satisfies the Boundary, as shown in Formula (19).

$$a_{min} \leq IVq - ROFHHM^{x,y}(a_1, a_2, \ldots a_n) \leq a_{max} \tag{19}$$

Proof: According to formula (6), the scoring function is defined as:

$$S(a_{\min}) = \frac{1}{2}[(\min_i\{u_{a_i}^-\})^q + (\min_i\{u_{a_i}^+\})^q - (\max_i\{v_{a_i}^-\})^q - (\max_i\{v_{a_i}^+\})^q]$$

$$S(a_{\max}) = \frac{1}{2}[(\max_i\{u_{a_i}^-\})^q + (\max_i\{u_{a_i}^+\})^q - (\min_i\{v_{a_i}^-\})^q - (\min_i\{v_{a_i}^+\})^q]$$

For any $u_{a_i}^-$、 $u_{a_i}^+$、 $v_{a_i}^-$、 $v_{a_i}^+$, meet:

$$(\min_i\{u_{a_i}^-\})^q \leq \left(u_{a_i}^-\right)^q \leq (\max_i\{u_{a_i}^-\})^q$$

$$(\min_i\{u_{a_i}^+\})^q \leq \left(u_{a_i}^+\right)^q \leq (\max_i\{u_{a_i}^+\})^q$$

$$(\min_i\{v_{a_i}^-\})^q \leq \left(v_{a_i}^-\right)^q \leq (\max_i\{v_{a_i}^-\})$$

$$(\min_i\{v_{a_i}^+\})^q \leq \left(v_{a_i}^+\right)^q \leq (\max_i\{v_{a_i}^+\})^q$$

conclude:

$$S(\min_i\{a_i\}) \leq S(IVq - ROFHHM) \leq S(\max_i\{a_i\})$$

According to the comparison and definition of interval values:

$$\min_i\{a_i\} \leq IVq - ROFHHM \leq \max_i\{a_i\}$$

Theorem 5: When $x, y, \varphi$ changes, the IVq-ROFHHM operator will also change, so three special cases of IVq-ROFHHM operator will be discussed below.

(1) When $\varphi = 1$, as shown in formula (20)

$$IVq - ROFHHM^{x,y}\left(a_1, a_2, \ldots a_n\right) = \frac{1}{x+y}\left(\underset{i=1,j=1}{\overset{n}{\otimes_H}}\left(xa_i \oplus_H ya_j\right)^{\frac{2}{n(n+1)}}\right)$$

$$= \left( \begin{array}{l} \left[ 1 - \left(1 - \prod_{i=1,j=1}^{n}\left(1 - \left(1 - \left(\mu_i^-\right)^q\right)^x \left(1 - \left(\mu_j^-\right)^q\right)^y\right)^{\frac{2}{n(n+2)}}\right)^{\frac{1}{x+y}}, 1 - \left(1 - \prod_{i=1,j=1}^{n}\left(1 - \left(1 - \left(\mu_i^+\right)^q\right)^x \left(1 - \left(\mu_j^+\right)^q\right)^y\right)^{\frac{2}{n(n+2)}}\right)^{\frac{1}{x+y}} \right], \\ \left[ \left(1 - \prod_{i=1,j=1}^{n}\left(1 - \left(\left(v_i^-\right)^q\right)^x \left(\left(v_j^-\right)^q\right)^y\right)^{\frac{2}{n(n+2)}}\right)^{\frac{1}{x+y}}, \left(1 - \prod_{i=1,j=1}^{n}\left(1 - \left(\left(v_i^-\right)^q\right)^x \left(\left(v_j^-\right)^q\right)^y\right)^{\frac{2}{n(n+2)}}\right)^{\frac{1}{x+y}} \right] \end{array} \right) \tag{20}$$

(2) When $x = 1, y = 0$, as shown in formula (21)

$$IVq - ROFHHM^{1,0}\left(a_1, a_2, \ldots a_n\right) = \frac{1}{x+y}\left(\underset{i=1,j=1}{\overset{n}{\otimes_H}}\left(xa_i \oplus_H ya_j\right)^{\frac{2}{n(n+1)}}\right)$$

$$\left(\left[\frac{\varphi\prod_{i=1,j=1}^{n}\left(V_{ij}^{+}-W_{ij}^{+}\right)^{\frac{2}{n(n+1)}}}{\varphi\prod_{i=1,j=1}^{n}\left(V_{ij}^{-}+\left(\varphi^{2}-1\right)W_{ij}^{-}\right)^{\frac{2}{n(n+1)}}+\left(\varphi^{2}-\varphi\right)\prod_{i=1,j=1}^{n}\left(V_{ij}^{-}-W_{ij}^{-}\right)^{\frac{2}{n(n+1)}}},\frac{\varphi\prod_{i=1,j=1}^{n}\left(V_{ij}^{+}-W_{ij}^{+}\right)^{\frac{2}{n(n+1)}}}{\varphi\prod_{i=1,j=1}^{n}\left(V_{ij}^{-}+\left(\varphi^{2}-1\right)W_{ij}^{-}\right)^{\frac{2}{n(n+1)}}+\left(\varphi^{2}-\varphi\right)\prod_{i=1,j=1}^{n}\left(V_{ij}^{-}-W_{ij}^{-}\right)^{\frac{2}{n(n+1)}}}\right],\left[\frac{\varphi\left(\prod_{i=1,j=1}^{n}\left(N_{ij}^{-}+\left(\varphi^{2}-1\right)M_{ij}^{-}\right)^{\frac{2}{n(n+1)}}-\prod_{i=1,j=1}^{n}\left(N_{ij}^{-}-M_{ij}^{-}\right)^{\frac{2}{n(n+1)}}\right)}{\varphi\prod_{i=1,j=1}^{n}\left(N_{ij}^{-}+\left(\varphi^{2}-1\right)M_{ij}^{-}\right)^{\frac{2}{n(n+1)}}+\left(\varphi^{2}+\varphi\right)\prod_{i=1,j=1}^{n}\left(N_{ij}^{-}-M_{ij}^{-}\right)^{\frac{2}{n(n+1)}}},\frac{\varphi\left(\prod_{i=1,j=1}^{n}\left(N_{ij}^{-}+\left(\varphi^{2}-1\right)M_{ij}^{-}\right)^{\frac{2}{n(n+1)}}-\prod_{i=1,j=1}^{n}\left(N_{ij}^{-}-M_{ij}^{-}\right)^{\frac{2}{n(n+1)}}\right)}{\varphi\prod_{i=1,j=1}^{n}\left(N_{ij}^{-}+\left(\varphi^{2}-1\right)M_{ij}^{-}\right)^{\frac{2}{n(n+1)}}+\left(\varphi^{2}+\varphi\right)\prod_{i=1,j=1}^{n}\left(N_{ij}^{-}-M_{ij}^{-}\right)^{\frac{2}{n(n+1)}}}\right]\right) \quad (21)$$

among

$$V_{\mathrm{ij}}=\varphi^{x}\left(1-(\varphi-1)\left(\mu_{i}^{-}\right)^{q}\right)^{x}$$

$$V_{ij}^{+}=\varphi^{x}\left(1-(\varphi-1)\left(\mu_{i}^{+}\right)^{q}\right)^{x}$$

$$W_{ij}^{-}=\left(\varphi\left(1-\left(\mu_{i}^{-}\right)^{q}\right)\right)^{x}$$

$$W_{ij}^{+}=\left(\varphi\left(1-\left(\mu_{i}^{+}\right)^{q}\right)\right)^{x}$$

$$N_{ij}^{-}=\left(1+\varphi^{2}+\varphi\left(1-\left(v_{i}^{-}\right)^{q}+\varphi\left(v_{i}^{-}\right)^{q}\right)\right)^{x}$$

$$N_{ij}^{+}=\left(1+\varphi^{2}+\varphi\left(1-\left(v_{i}^{+}\right)^{q}+\varphi\left(v_{i}^{+}\right)^{q}\right)\right)^{x}$$

$$M_{ij}^{-}=\left(\varphi\left(v_{i}^{-}\right)^{q}\right)^{x}$$

$$M_{ij}^{+}=\left(\varphi\left(v_{i}^{+}\right)^{q}\right)^{x}$$

(3) When $x=0, y=1$, as shown in formula (22)

$$IVq-ROFHHM^{1,0}\left(a_{1},a_{2},\ldots a_{n}\right)=\frac{1}{x+y}\left(\underset{i=1,j=1}{\overset{n}{\otimes_{H}}}\left(xa_{i}\oplus_{H} ya_{j}\right)^{\frac{2}{n(n+1)}}\right)$$

$$\left(\left[\frac{\varphi\prod_{i=1,j=1}^{n}\left(V_{ij}^{+}-W_{ij}^{+}\right)^{\frac{2}{n(n+1)}}}{\varphi\prod_{i=1,j=1}^{n}\left(V_{ij}^{-}+\left(\varphi^{2}-1\right)W_{ij}^{-}\right)^{\frac{2}{n(n+1)}}+\left(\varphi^{2}-\varphi\right)\prod_{i=1,j=1}^{n}\left(V_{ij}^{-}-W_{ij}^{-}\right)^{\frac{2}{n(n+1)}}},\frac{\varphi\prod_{i=1,j=1}^{n}\left(V_{ij}^{+}-W_{ij}^{+}\right)^{\frac{2}{n(n+1)}}}{\varphi\prod_{i=1,j=1}^{n}\left(V_{ij}^{-}+\left(\varphi^{2}-1\right)W_{ij}^{-}\right)^{\frac{2}{n(n+1)}}+\left(\varphi^{2}-\varphi\right)\prod_{i=1,j=1}^{n}\left(V_{ij}^{-}-W_{ij}^{-}\right)^{\frac{2}{n(n+1)}}}\right],\left[\frac{\varphi\left(\prod_{i=1,j=1}^{n}\left(N_{ij}^{-}+\left(\varphi^{2}-1\right)M_{ij}^{-}\right)^{\frac{2}{n(n+1)}}-\prod_{i=1,j=1}^{n}\left(N_{ij}^{-}-M_{ij}^{-}\right)^{\frac{2}{n(n+1)}}\right)}{\varphi\prod_{i=1,j=1}^{n}\left(N_{ij}^{-}+\left(\varphi^{2}-1\right)M_{ij}^{-}\right)^{\frac{2}{n(n+1)}}+\left(\varphi^{2}+\varphi\right)\prod_{i=1,j=1}^{n}\left(N_{ij}^{-}-M_{ij}^{-}\right)^{\frac{2}{n(n+1)}}},\frac{\varphi\left(\prod_{i=1,j=1}^{n}\left(N_{ij}^{-}+\left(\varphi^{2}-1\right)M_{ij}^{-}\right)^{\frac{2}{n(n+1)}}-\prod_{i=1,j=1}^{n}\left(N_{ij}^{-}-M_{ij}^{-}\right)^{\frac{2}{n(n+1)}}\right)}{\varphi\prod_{i=1,j=1}^{n}\left(N_{ij}^{-}+\left(\varphi^{2}-1\right)M_{ij}^{-}\right)^{\frac{2}{n(n+1)}}+\left(\varphi^{2}+\varphi\right)\prod_{i=1,j=1}^{n}\left(N_{ij}^{-}-M_{ij}^{-}\right)^{\frac{2}{n(n+1)}}}\right]\right) \quad (22)$$

among

$$V_{ij}^{-} = \varphi^{y}\left(1-(\varphi-1)\left(\mu_i^{-}\right)^{q}\right)^{y}$$

$$V_{ij}^{+} = \varphi^{y}\left(1-(\varphi-1)\left(\mu_i^{+}\right)^{q}\right)^{y}$$

$$W_{ij}^{-} = \left(\varphi\left(1-\left(\mu_i^{-}\right)^{q}\right)\right)^{y}$$

$$W_{ij}^{+} = \left(\varphi\left(1-\left(\mu_i^{+}\right)^{q}\right)\right)^{y}$$

$$N_{ij}^{-} = \left(1+\varphi^{2}+\varphi\left(1-\left(v_i^{-}\right)^{q}+\varphi\left(v_i^{-}\right)^{q}\right)\right)^{y}$$

$$N_{ij}^{+} = \left(1+\varphi^{2}+\varphi\left(1-\left(v_i^{+}\right)^{q}+\varphi\left(v_i^{+}\right)^{q}\right)\right)^{y}$$

$$M_{ij}^{-} = \left(\varphi\left(v_i^{-}\right)^{q}\right)^{y}$$

$$M_{ij}^{+} = \left(\varphi\left(v_i^{+}\right)^{q}\right)^{y}$$

# 4. Interval-valued generalized orthogonal HHM operator

According to the definition of Hamacher operation in Section 2.2, this section proposes an interval valued q-rung Orthopair Fuzzy Hamacher operation, which is used to aggregate the decision matrix of each expert. The algorithm of IVq-ROFHHM operator is given in section 4.1, and the definition and properties of IVq-ROFHHMWA operator are expounded in section 4.2.

## 4.1 IVq-ROFHHMWA operator

Definition: Let $x, y \geq 0$ which are not all 0 at the same time, $a_i = \langle [u_i^-, u_i^+], [v_i^-, v_i^+] \rangle (i = 1,2, \dots n)$ is a group of IVq-ROFNs, $w = (w_1, w_1, \dots, w_n)^T$ are weight vectors, $w_i \in [0,1], \sum_{i=1}^{n} w_i = 1$. There is a mapping: $IVq - ROFHHMWA^n \to IVq - ROFHHMWA$, then the IVq-ROFWHHM operator is shown in Formula (23).

$$IVq - ROFHHMWA_w^{x,y}\left(a_1, a_2, \dots a_n\right) = \frac{1}{x+y}\left(\underset{i=1,j=1}{\overset{n}{\otimes_H}}\left(x a_i^{w_i} \oplus_H y a_j^{w_j}\right)^{\frac{2}{n(n+1)}}\right) \quad (23)$$

Theorem 6: Let $x, y \geq 0$, and not be 0 at the same time, $a_i = \langle [u_i^-, u_i^+], [v_i^-, v_i^+] \rangle (i = 1,2, \dots n)$ is a group of IVq-ROFNs, then the operation formula of IVq-ROFWHHM operator is shown in (24).

$$IVq - ROFHHMWA_w^{x,y}\left(a_1, a_2, \dots a_n\right) = \frac{1}{x+y}\left(\underset{i=1,j=1}{\overset{n}{\otimes_H}}\left(x a_i^{w_i} \oplus_H y a_j^{w_j}\right)^{\frac{2}{n(n+1)}}\right)$$

$$=\left(\begin{array}{l}\left[\sqrt[q]{\dfrac{a^- - b^-}{a^- + (\varphi-1)b^-}}, \sqrt[q]{\dfrac{a^+ - b^+}{a^+ + (\varphi-1)b^+}}\right], \\ \left[\sqrt[q]{\dfrac{\varphi\left(d^- - c^-\right)^{\frac{1}{x+y}}}{\left(c^- + (\varphi^2-1)d^-\right)^{\frac{1}{x+y}} + (\varphi-1)\left(d^- - c^-\right)^{\frac{1}{x+y}}}}, \sqrt[q]{\dfrac{\varphi\left(d^+ - c^+\right)^{\frac{1}{x+y}}}{\left(c^+ + (\varphi^2-1)d^+\right)^{\frac{1}{x+y}} + (\varphi-1)\left(d^+ - c^+\right)^{\frac{1}{x+y}}}}\right]\end{array}\right) \tag{24}$$

In formula (24), some formulas are expressed by formulas (25)-(32).

$$a^- = \left(\prod_{i=1,j=1}^{n}\left(V_{ij}^- + (\varphi^2-1)W_{ij}^-\right)^{\frac{2}{n(n+1)}} + (\varphi^2-1)\prod_{i=1,j=1}^{n}\left(V_{ij}^- - W_{ij}^-\right)^{\frac{2}{n(n+1)}}\right)^{\frac{1}{x+y}} \tag{25}$$

$$a^+ = \left(\prod_{i=1,j=1}^{n}\left(V_{ij}^+ + (\varphi^2-1)W_{ij}^+\right)^{\frac{2}{n(n+1)}} + (\varphi^2-1)\prod_{i=1,j=1}^{n}\left(V_{ij}^+ - W_{ij}^+\right)^{\frac{2}{n(n+1)}}\right)^{\frac{1}{x+y}} \tag{26}$$

$$b^- = \left(\prod_{i=1,j=1}^{n}\left(V_{ij}^- + (\varphi^2-1)W_{ij}^-\right)^{\frac{2}{n(n+1)}} - \prod_{i=1,j=1}^{n}\left(V_{ij}^- - W_{ij}^-\right)^{\frac{2}{n(n+1)}}\right)^{\frac{1}{x+y}} \tag{27}$$

$$b^+ = \left(\prod_{i=1,j=1}^{n}\left(V_{ij}^+ + (\varphi^2-1)W_{ij}^+\right)^{\frac{2}{n(n+1)}} - \prod_{i=1,j=1}^{n}\left(V_{ij}^+ - W_{ij}^+\right)^{\frac{2}{n(n+1)}}\right)^{\frac{1}{x+y}} \tag{28}$$

$$c^- = \left(\prod_{i=1,j=1}^{n}\left(N_{ij}^- + (\varphi^2-1)M_{ij}^-\right)^{\frac{2}{n(n+1)}} - \prod_{i=1,j=1}^{n}\left(N_{ij}^- - M_{ij}^-\right)^{\frac{2}{n(n+1)}}\right)^{\frac{1}{x+y}} \tag{29}$$

$$c^+ = \left(\prod_{i=1,j=1}^{n}\left(N_{ij}^+ + (\varphi^2-1)M_{ij}^+\right)^{\frac{2}{n(n+1)}} - \prod_{i=1,j=1}^{n}\left(N_{ij}^+ - M_{ij}^+\right)^{\frac{2}{n(n+1)}}\right)^{\frac{1}{x+y}} \tag{30}$$

$$\begin{aligned} d^- = &\left(\prod_{i=1,j=1}^{n}\left(N_{ij}^- + (\varphi^2-1)M_{ij}^-\right)^{\frac{2}{n(n+1)}} + (\varphi^2-1)\prod_{i=1,j=1}^{n}\left(N_{ij}^- - M_{ij}^-\right)^{\frac{2}{n(n+1)}}\right)^{\frac{1}{x+y}} \\ &+ (\varphi-1)\left(\prod_{i=1,j=1}^{n}\left(N_{ij}^- + (\varphi^2-1)M_{ij}^-\right)^{\frac{2}{n(n+1)}} + \prod_{i=1,j=1}^{n}\left(N_{ij}^- - M_{ij}^-\right)^{\frac{2}{n(n+1)}}\right)^{\frac{1}{x+y}} \end{aligned} \tag{31}$$

$$\begin{aligned} d^+ = &\left(\prod_{i=1,j=1}^{n}\left(N_{ij}^+ + (\varphi^2-1)M_{ij}^+\right)^{\frac{2}{n(n+1)}} + (\varphi^2-1)\prod_{i=1,j=1}^{n}\left(N_{ij}^+ - M_{ij}^+\right)^{\frac{2}{n(n+1)}}\right)^{\frac{1}{x+y}} \\ &+ (\varphi-1)\left(\prod_{i=1,j=1}^{n}\left(N_{ij}^+ + (\varphi^2-1)M_{ij}^+\right)^{\frac{2}{n(n+1)}} + \prod_{i=1,j=1}^{n}\left(N_{ij}^+ - M_{ij}^+\right)^{\frac{2}{n(n+1)}}\right)^{\frac{1}{x+y}} \end{aligned} \tag{32}$$

among

$$V_{ij}^{-}=\left(\left(1+(\varphi-1)\left(1-\left(\mu_i^{-}\right)^q\right)\right)^{w_i}+\left(\varphi^2-1\right)\left(\left(\mu_i^{-}\right)^q\right)^{w_i}\right)^x\left(\left(1+(\varphi-1)\left(1-\left(\mu_j^{-}\right)^q\right)\right)^{w_j}+\left(\varphi^2-1\right)\left(\left(\mu_j^{-}\right)^q\right)^{w_j}\right)^y$$

$$V_{ij}^{+}=\left(\left(1+(\varphi-1)\left(1-\left(\mu_i^{+}\right)^q\right)\right)^{w_i}+\left(\varphi^2-1\right)\left(\left(\mu_i^{+}\right)^q\right)^{w_i}\right)^x\left(\left(1+(\varphi-1)\left(1-\left(\mu_j^{+}\right)^q\right)\right)^{w_j}+\left(\varphi^2-1\right)\left(\left(\mu_j^{+}\right)^q\right)^{w_j}\right)^y$$

$$W_{ij}^{-}=\left(\left(1+(\varphi-1)\left(1-\left(\mu_i^{-}\right)^q\right)\right)^{w_i}-\left(\left(\mu_i^{-}\right)^q\right)^{w_j}\right)^x\left(\left(1+(\varphi-1)\left(1-\left(\mu_j^{-}\right)^q\right)\right)^{w_j}-\left(\left(\mu_j^{-}\right)^q\right)^{w_j}\right)^y$$

$$W_{ij}^{+}=\left(\left(1+(\varphi-1)\left(1-\left(\mu_i^{+}\right)^q\right)\right)^{w_i}-\left(\left(\mu_i^{+}\right)^q\right)^{w_j}\right)^x\left(\left(1+(\varphi-1)\left(1-\left(\mu_j^{+}\right)^q\right)\right)^{w_j}-\left(\left(\mu_j^{+}\right)^q\right)^{w_j}\right)^y$$

$$N_{ij}^{-}=\left(\left(1+(\varphi-1)\left(1-\left(v_i^{-}\right)^q\right)\right)^{w_i}+\left(\varphi^2-1\right)\left(1-\left(v_i^{-}\right)^q\right)^{w_i}\right)^x\left(\left(1+(\varphi-1)\left(v_j^{-}\right)^q\right)^{w_j}+\left(\varphi^2-1\right)\left(1-\left(v_j^{-}\right)^q\right)^{wj}\right)^y$$

$$N_{ij}^{+}=\left(\left(1+(\varphi-1)\left(1-\left(v_i^{+}\right)^q\right)\right)^{w_i}+\left(\varphi^2-1\right)\left(1-\left(v_i^{+}\right)^q\right)^{w_i}\right)^x\left(\left(1+(\varphi-1)\left(v_j^{+}\right)^q\right)^{w_j}+\left(\varphi^2-1\right)\left(1-\left(v_j^{+}\right)^q\right)^{wj}\right)^y$$

$$M_{ij}^{-}=\left(\left(1+(\varphi-1)\left(v_i^{-}\right)^q\right)^{w_i}-\left(1-\left(v_i^{-}\right)^q\right)^{w_i}\right)^x\left(\left(1+(\varphi-1)\left(v_j^{-}\right)^q\right)^{w_j}-\left(1-\left(v_j^{-}\right)^q\right)^{wj}\right)^y$$

$$M_{ij}^{+}=\left(\left(1+(\varphi-1)\left(v_i^{+}\right)^q\right)^{w_i}-\left(1-\left(v_i^{+}\right)^q\right)^{w_i}\right)^x\left(\left(1+(\varphi-1)\left(v_j^{+}\right)^q\right)^{w_j}-\left(1-\left(v_j^{+}\right)^q\right)^{wj}\right)^y$$

According to the derivation and proof of formula (16), theorem 6 can also be proved, and this paper will not give a detailed proof.

IVq-ROFHHMWA operator satisfies permutation, boundary,monotonicity. Below, we give the properties of IVq-ROFHHMWA operator.

Theorem 7 (Permutation): Let $(\beta_1,\beta_2,\dots\beta_n)$ be any permutation of $(a_1,a_2,\dots a_n)$, then IVq-ROFHHMWA operator satisfies permutation, as shown in formula (33).

$$IVq-ROFHHMWA^{x,y}(\beta_1,\beta_2,\dots\beta_n)=IVq-ROFWHHMWA^{x,y}(a_1,a_2,\dots a_n) \tag{33}$$

Theorem 8 (Monotonicity): Let $a_i=\langle[u_i^-,u_i^+],[v_i^-,v_i^+]\rangle$ and $a_i'=\langle[u_i^{-\prime},u_i^{+\prime}],[v_i^{-\prime},v_i^{+\prime}]\rangle(i=1,2,\dots,\mathrm{n})$ are two groups of IVq-ROFNs, and satisfy: $u_i^-\le u_i^{-\prime}$, $u_i^+\le u_i^{+\prime}$, $v_i^-\ge v_i^{-\prime}$, $v_i^+\ge v_i^{+\prime}$, then the IVq-ROFEHM operator satisfies monotonicity, as shown in formula (34).

$$IVq-ROFHHMWA^{x,y}(a_1,a_2,\dots a_n)\le IVq-ROFWHHMWA^{x,y}(a_1',a_2'\dots a_n') \tag{34}$$

Theorem 9 (Boundary): Let $(a_1,a_2,\dots a_n)$ be a set of IVq-ROFNs, $a_{\min}=\min_i\{a_i\}$, $a_{\max}=\max_i\{a_i\}$. Then the IVq-ROFWEHM operator satisfies the Boundary, as shown in Formula (35).

$$a_{min}\le IVq-ROFHHMWA^{x,y}(a_1,a_2,\dots a_n)\le a_{max} \tag{35}$$

The proof process of formulas (33) to (35) can be easily deduced and proved by using formulas (17), (18) and (19), so the detailed proof process will not be given in this paper.

### 4.2 IVq-ROFHHMGA operator

Definition 10: Let x, y ≥ 0, and not equal to 0 at the same time, $a_i = \langle [u_i^-, u_i^+], [v_i^-, v_i^+] \rangle (i = 1,2, \ldots n)$ is a group of IVq-ROFNs, and $w = (w_1, w_1, \ldots, w_n)^T$ is the weight vector, $w_i \in [0,1]$, $\sum_{i=1}^{n} w_i = 1$. There is a mapping: $IVq - ROFHHMGA^n \to IVq - ROFHHMGA$, then the IVq-ROFHHMGA operator is shown in Formula (36).

$$IVq - ROFHHMGA_w^{x,y}(a_1, a_2, \ldots a_n) = \frac{1}{x+y}\left( \underset{i=1,j=1}{\overset{n}{\oplus_H}} \left( x a_i^{w_i} \otimes_H y a_j^{w_j} \right)^{\frac{2}{n(n+1)}} \right) \tag{36}$$

Theorem 10: Let x, y≥0, and not equal to 0 at the same time, $a_i = \langle [u_i^-, u_i^+], [v_i^-, v_i^+] \rangle (i = 1,2, \ldots n)$ is a group of IVq-ROFNs, and the operation formula of IVq-ROFWHHM operator is shown in (37).

$$IVq - ROFHHMGA_w^{x,y}(a_1, a_2, \ldots a_n) = \frac{1}{x+y}\left( \underset{i=1,j=1}{\overset{n}{\oplus_H}} \left( x a_i^{w_i} \otimes_H y a_j^{w_j} \right)^{\frac{2}{n(n+1)}} \right)$$

$$= \left( \begin{array}{l} \left[ \sqrt[q]{\frac{\varphi\left(a^- - b^-\right)^{\frac{1}{x+y}}}{\left(a^- + \left(\varphi^2 - 1\right)b^-\right)^{\frac{1}{x+y}} + (\varphi - 1)\left(a^- - b^-\right)^{\frac{1}{x+y}}}}, \sqrt[q]{\frac{\varphi\left(a^+ - b^+\right)^{\frac{1}{x+y}}}{\left(a^+ + \left(\varphi^2 - 1\right)b^+\right)^{\frac{1}{x+y}} + (\varphi - 1)\left(a^+ - b^+\right)^{\frac{1}{x+y}}}} \right], \\ \left[ \sqrt[q]{\frac{d^- - c^-}{c^- + (\varphi - 1)d^-}}, \sqrt[q]{\frac{d^+ - c^+}{c^+ + (\varphi - 1)d^+}} \right] \end{array} \right) \tag{37}$$

In formula (37), some formulas are expressed by formulas (25)-(32).

According to the derivation and proof of formula (16), theorem 6 can also be proved, and this paper will not give a detailed proof.

IVq-ROFHHMGA operator satisfies permutation, boundary,, monotonicity. Below, we give the properties of IVq-ROFHHMGA operator.

Theorem 11 (Permutation): If $(\beta_1, \beta_2, \ldots \beta_n)$ is any permutation of $(a_1, a_2, \ldots a_n)$, then IVq-ROFHHMGA operator satisfies permutation, as shown in formula (38).

$$IVq - ROFHHMGA^{x,y}(\beta_1, \beta_2, \ldots \beta_n) = IVq - ROFWHHMGA^{x,y}(a_1, a_2, \ldots a_n) \tag{38}$$

Theorem 12 (Monotonicity): Let $a_i = \langle [u_i^-, u_i^+], [v_i^-, v_i^+] \rangle$ and $a_i' = \langle [u_i^{-\prime}, u_i^{+\prime}], [v_i^{-\prime}, v_i^{+\prime}] \rangle (i = 1,2, \ldots, n)$ be two groups of IVq-ROFNs, and satisfy: $u_i^- \le u_i^{-\prime}$, $u_i^+ \le u_i^{+\prime}$, $v_i^- \ge v_i^{-\prime}$, $v_i^+ \ge v_i^{+\prime}$, then the IVq-ROFEHHMGA operator satisfies monotonicity, as shown in formula (39).

$$IVq - ROFHHMGA^{x,y}(a_1, a_2, \ldots a_n) \le IVq - ROFWHHMGA^{x,y}(a_1', a_2' \ldots a_n') \tag{39}$$

Theorem 13 (Boundary): Let $(a_1, a_2, \ldots a_n)$ be a set of IVq-ROFNs, $a_{\min} = \min_i \{a_i\}$, $a_{\max} = \max_i \{a_i\}$. Then the IVq-ROFHHMGA operator satisfies the Boundary, as shown in Formula (40).

$$a_{min} \le IVq - ROFHHMGA^{x,y}(a_1, a_2, \ldots a_n) \le a_{max} \tag{40}$$

The proving process of formulas (38) to (40) can be easily deduced and proved by using formulas (17), (18) and (19), so the detailed proving process will not be given in this paper.

## 5. Group decision-making method

### 5.1 generalized orthogonal group decision-making conditions

In order to select the best scheme from the alternatives or sort the multiple schemes more scientifically, it is often necessary to invite more experts to participate in the decision-making. In order to improve the efficiency of decision-making and reduce the cost given by the decision matrix in the decision-making process, the interval valued q-rung Orthopair Fuzzy number proposed in this paper gives the experts more degrees of freedom, and does not need to consider whether the weights of experts and the attributes are given. Therefore, the environmental conditions of group decision-making in this paper are described in detail as follows.

In daily life, it is a common decision-making problem to choose the best scheme from M schemes with N attributes with T experts participating, and the weights of experts and attributes are unknown. Its mathematical expression is: the set of K experts is: $E = \{e_1, e_2, \ldots, e_k\}$, and the attribute weights satisfy: $\sum_{i=1}^{t} \lambda_i = 1$， $\lambda_i \geq 0, i \in \{1,2,\ldots,t\}$; M alternative schemes are: $X = \{x_1, x_2, \ldots, x_m\}$, and the attribute set of each scheme is: $C = \{c_1, c_2, \ldots, c_n\}$, with attribute weight $\sum_{i=1}^{n} w_i = 1, w_i \geq 0, i \in \{1,2,\ldots,n\}$. In order to make experts give decision matrix more flexibly, all households are required to give interval valued q-rung Orthopair Fuzzy decision matrix, in which the interval valued q-rung Orthopair Fuzzy decision matrix given by the t expert is: $A^{(t)} = (a_{ij}^{(t)})_{m\times n}$, and the decision matrix is shown in formula (49).

$$A^{(t)} = \begin{bmatrix} a_{11}^{(t)} \cdots a_{1j}^{(t)} \cdots a_{1n}^{(t)} \\ \ddots \quad \ddots \ \vdots \\ a_{mj}^{(t)} \cdots a_{ij}^{(t)} \cdots a_{in}^{(t)} \\ \vdots \ \ddots \ \vdots \ \ddots \ \vdots \\ a_{m1}^{(t)} \cdots a_{mj}^{(t)} \cdots a_{mn}^{(t)} \end{bmatrix}_{m\times n}, i = 1,2,\cdots,m, j = 1,2,\cdots,n \quad (41)$$

Among them, the interval-valued generalized orthogonal fuzzy numbers $a_{ij}^{(t)} =<$ $[u^{(t)}{}_{\mathrm{ij}}^{-}, u^{(t)}{}_{\mathrm{ij}}^{+}], [v^{(t)}{}_{\mathrm{ij}}^{-}, v^{(t)}{}_{\mathrm{ij}}^{+}] >$ can find appropriate value of q, so that all fuzzy numbers in the decision matrix given by experts satisfy: $0 \leq (u^{(t)}{}_{ij}^{-})^q \leq (u^{(t)}{}_{ij}^{+})^q \leq 1, 0 \leq (v^{(t)}{}_{ij}^{-})^q \leq (v^{(t)}{}_{ij}^{+})^q \leq 1, 0 \leq (u^{(t)}{}_{ij}^{-})^q + (v^{(t)}{}_{ij}^{-})^q \leq (u^{(t)}{}_{ij}^{+})^q + (v^{(t)}{}_{ij}^{+})^q \leq 1$.

In order to make it convenient for t experts to get the order of schemes and the best scheme efficiently and accurately, the decision-making method needs to use IVq-ROFHHMWA proposed in the fourth part of this paper to assemble the expert decision matrix; The score is calculated by using the score function given by formula (6), and the ranking and the best scheme can be obtained according to the score.

## 5.2 Group decision-making method

According to the decision-making environment described in section 5.1, experts only need to care about the decision-making judgment of alternatives, without considering whether the decision value given by themselves is related to other experts or the relationship between the attributes of alternatives. They only need to give their own interval valued q-rung Orthopair Fuzzy decision matrix, thus greatly reducing the decision-making cost. In this decision-making environment, firstly, the Q value satisfying the interval valued q-rung Orthopair Fuzzy number should be obtained according to the matrix given by the decision maker, and the follow-up work can be carried out under the condition of determining Q. According to the decision matrix given by experts, the attitude of experts can be reflected laterally, and reasonable expert weights can be set. According to expert weights and decision matrix, IVq-ROFHHMWA operator can be used to obtain aggregation matrix. At this time, the score function proposed in this paper is used to calculate the score of each scheme, and the ranking of each scheme and the best scheme can be obtained according to the score. For experts, this decision-making algorithm only cares about the problems and alternatives, while for decision-making units, it only needs to consider the problems and the experts needed, thus reducing the decision-making cost. The decision-making algorithm proposed in this paper includes six steps.

Step 1: Determine the value of Q, traverse the elements in each expert decision matrix in turn, and find out the value that satisfies the generalized orthogonal condition of interval value Q, that is, the value of Q should satisfy: $\left(u^{+}_{a_{ij}^{(k)}}\right)^q + \left(v^{+}_{a_{ij}^{(k)}}\right)^q \le 1, \left(u^{-}_{a_{ij}^{(k)}}\right)^q + \left(v^{-}_{a_{ij}^{(k)}}\right)^q \le 1$.

Step 2: Give appropriate expert weights. $\lambda = \left(\lambda_1, \lambda_2, ..., \lambda_n\right)^T$ 。

Step 3: according to the weight of the experts in step 2, the decision matrix of different experts is assembled by using formula (42) to obtain the interval valued q-rung Orthopair Fuzzy number matrix. $\boldsymbol{R} = (\boldsymbol{a_{ij}})_{\boldsymbol{m\times n}}$

$$r_{ij} = \mathrm{IVq-ROFHHMWA}\left(a_{ij}^{(1)}, a_{ij}^{(2)}, \cdots, a_{ij}^{(t)}\right) \quad i = 1, 2, \cdots m\ ;\ \ j = 1, 2, \cdots n \qquad (42)$$

Step 4: According to the given attribute weight $\omega = \left(\omega_1, \omega_2, ..., \omega_n\right)^T$, assemble the R scheme by using formula (43) to obtain the scheme $x_i$.

$$x_i = \mathrm{IVq-ROFHHMWA}\left(r_{i1}, r_{i2}, \cdots, r_{in}\right) \qquad (43)$$

Step 5: Use the score function (9) and the exact value function (10) to calculate the $x_i$, and sort them according to their size.

Step 6: According to the sorting result of $x_i$, the maximum value is the optimal scheme.

# 6.Multiple attributes decision making with the proposed operators and illustrated examples

## 6.1 Number Case

For experts' more accurate judgment, interval valued q-rung Orthopair Fuzzy numbers are used to express the scores of each index. Experts suggest that five indexes be adopted, namely: C1, C2, C3, C4, C5. Three experts give decision matrix, as shown Table1 to Table3.

Table 1 Decision matrix given by expert 1$E_1$

| | C1 (genetic) | C2 (blood pressure measurement) | C3 (Eating habits) | C4 (Exercise Habit) | C5 (stress of life) |
|---|---|---|---|---|---|
| A1 | ([0.35, 0.45], [0.5, 0.65]) | ([0.8, 0.85], [0.15, 0.2]) | ([0.6, 0.7], [0.3, 0.4]) | ([0.7, 0.8], [0.2, 0.3]) | ([0.65, 0.7], [0.35, 0.4]) |
| A2 | ([0.55, 0.6], [0.4, 0.5]) | ([0.9, 0.95], [0.1, 0.2]) | ([0.75, 0.85], [0.2, 0.3]) | ([0.85, 0.9], [0.1, 0.15]) | ([0.75, 0.8], [0.2, 0.3]) |
| A3 | ([0.4, 0.5], [0.5, 0.6]) | ([0.85, 0.9], [0.1, 0.2]) | ([0.75, 0.8], [0.2, 0.3]) | ([0.8, 0.9], [0.1, 0.2]) | ([0.65, 0.75], [0.35, 0.4]) |
| A4 | ([0.35, 0.4], [0.6, 0.65]) | ([0.75, 0.85], [0.25, 0.3]) | ([0.6, 0.65], [0.25, 0.3]) | ([0.7, 0.8], [0.2, 0.3]) | ([0.55, 0.65], [0.3, 0.4]) |
| A5 | ([0.1, 0.2], [0.85, 0.9]) | ([0.65, 0.75], [0.3, 0.45]) | ([0.55, 0.6], [0.4, 0.5]) | ([0.6, 0.7], [0.3, 0.4]) | ([0.5, 0.6], [0.4, 0.5]) |

Table 2 Decision matrix given by expert 2$E_2$

| | C1 (genetic) | C2 (blood pressure measurement) | C3 (Eating habits) | C4 (Exercise Habit) | C5 (stress of life) |
|---|---|---|---|---|---|
| A1 | ([0.4, 0.45], [0.5, 0.6]) | ([0.85, 0.9], [0.1, 0.2]) | ([0.65, 0.75], [0.3, 0.4]) | ([0.75, 0.8], [0.2, 0.3]) | ([0.6, 0.7], [0.2, 0.3]) |
| A2 | ([0.5, 0.6], [0.4, 0.5]) | ([0.9, 0.95], [0.15, 0.2]) | ([0.75, 0.8], [0.2, 0.3]) | ([0.85, 0.95], [0.1, 0.15]) | ([0.7, 0.8], [0.2, 0.3]) |
| A3 | ([0.4, 0.5], [0.5, 0.6]) | ([0.85, 0.95], [0.1, 0.2]) | ([0.7, 0.8], [0.2, 0.3]) | ([0.85, 0.9], [0.1, 0.2]) | ([0.6, 0.75], [0.3, 0.35]) |
| A4 | ([0.35, 0.45], [0.55, 0.65]) | ([0.7, 0.85], [0.2, 0.3]) | ([0.65, 0.7], [0.3, 0.45]) | ([0.7, 0.8], [0.2, 0.3]) | ([0.6, 0.65], [0.3, 0.4]) |
| A5 | ([0.3, 0.4], [0.6, 0.7]) | ([0.7, 0.85], [0.2, 0.3]) | ([0.55, 0.65], [0.35, 0.4]) | ([0.6, 0.7], [0.3, 0.35]) | ([0.5, 0.6], [0.4, 0.5]) |

Table 3 Decision matrix given by expert 3$E_3$

| | C1 (genetic) | C2 (blood pressure measurement) | C3 (Eating habits) | C4 (Exercise Habit) | C5 (stress of life) |
|---|---|---|---|---|---|
| A1 | ([0.4, 0.5], [0.5, 0.6]) | ([0.8, 0.85], [0.15, 0.25]) | ([0.7, 0.8], [0.25, 0.3]) | ([0.75, 0.85], [0.2, 0.3]) | ([0.6, 0.7], [0.35, 0.4]) |
| A2 | ([0.5, 0.6], [0.4, 0.5]) | ([0.9, 0.95], [0.1, 0.2]) | ([0.75, 0.85], [0.2, 0.35]) | ([0.85, 0.95], [0.1, 0.15]) | ([0.7, 0.8], [0.2, 0.3]) |
| A3 | ([0.5, 0.65], [0.4, 0.5]) | ([0.85, 0.95], [0.1, 0.2]) | ([0.75, 0.8], [0.2, 0.3]) | ([0.85, 0.9], [0.1, 0.15]) | ([0.7, 0.8], [0.2, 0.3]) |
| A4 | ([0.35, 0.4], [0.6, 0.65]) | ([0.75, 0.85], [0.15, 0.25]) | ([0.6, 0.7], [0.3, 0.4]) | ([0.7, 0.8], [0.2, 0.3]) | ([0.6, 0.65], [0.3, 0.35]) |
| A5 | ([0.3, 0.4], [0.5, 0.6]) | ([0.7, 0.8], [0.2, 0.3]) | ([0.5, 0.65], [0.4, 0.5]) | ([0.65, 0.7], [0.2, 0.3]) | ([0.5, 0.6], [0.4, 0.5]) |

(1) Traverse the generalized orthogonal fuzzy numbers in all experts' evidence, and you can get q=3.

(2) Take $x = 3, y = 3, \varphi = 3$ and reasonably give the weight of each expert with a small difference, $\lambda = (0.330\ 0.334\ 0.336)$.

(3) Use formula (42) to aggregate the expert decisions $E^{(t)}$, and get R.

$$R = \begin{bmatrix} \begin{pmatrix}[0.79,0.83],\\ [0.96,0.96]\end{pmatrix}, & \begin{pmatrix}[0.94,0.96],\\ [0.96,0.96]\end{pmatrix}, & \begin{pmatrix}[0.90,0.92],\\ [0.96,0.96]\end{pmatrix}, & \begin{pmatrix}[0.92,0.94],\\ [0.96,0.96]\end{pmatrix}, & \begin{pmatrix}[0.89,0.91],\\ [0.96,0.96]\end{pmatrix} \\ \begin{pmatrix}[0.85,0.88],\\ [0.96,0.96]\end{pmatrix}, & \begin{pmatrix}[0.97,0.98],\\ [0.96,0.96]\end{pmatrix}, & \begin{pmatrix}[0.93,0.95],\\ [0.96,0.96]\end{pmatrix}, & \begin{pmatrix}[0.95,0.98],\\ [0.96,0.96]\end{pmatrix}, & \begin{pmatrix}[0.92,0.94],\\ [0.96,0.96]\end{pmatrix} \\ \begin{pmatrix}[0.82,0.86],\\ [0.96,0.96]\end{pmatrix}, & \begin{pmatrix}[0.95,0.98],\\ [0.96,0.96]\end{pmatrix}, & \begin{pmatrix}[0.92,0.94],\\ [0.96,0.96]\end{pmatrix}, & \begin{pmatrix}[0.95,0.97],\\ [0.96,0.96]\end{pmatrix}, & \begin{pmatrix}[0.89,0.93],\\ [0.96,0.96]\end{pmatrix} \\ \begin{pmatrix}[0.78,0.81],\\ [0.96,0.96]\end{pmatrix}, & \begin{pmatrix}[0.92,0.95],\\ [0.96,0.96]\end{pmatrix}, & \begin{pmatrix}[0.89,0.91],\\ [0.96,0.96]\end{pmatrix}, & \begin{pmatrix}[0.91,0.94],\\ [0.96,0.96]\end{pmatrix}, & \begin{pmatrix}[0.88,0.90],\\ [0.96,0.96]\end{pmatrix} \\ \begin{pmatrix}[0.65,0.74],\\ [0.96,0.96]\end{pmatrix}, & \begin{pmatrix}[0.91,0.94],\\ [0.96,0.96]\end{pmatrix} & \begin{pmatrix}[0.86,0.89],\\ [0.96,0.96]\end{pmatrix}, & \begin{pmatrix}[0.89,0.91],\\ [0.96,0.96]\end{pmatrix}, & \begin{pmatrix}[0.85,0.88],\\ [0.96,0.96]\end{pmatrix} \end{bmatrix}$$

(4) Based on the case attribute weight $\omega = \left(0.194,\ 0.234\ 0.218,\ 0.183,\ 0.171\right)^T$, the R scheme is aggregated by formula (43), and the interval valued q-rung Orthopair Fuzzy number of the scheme is obtained.,

$$x_1 = ([0.98, 0.98], [0.96, 0.96])$$
$$x_2 = ([0.98, 0.99], [0.96, 0.96])$$
$$x_3 = ([0.98, 0.98], [0.96, 0.96])$$
$$x_4 = ([0.97, 0.98], [0.96, 0.96])$$
$$x_5 = ([0.96, 0.97], [0.96, 0.96])$$

(5) Use formula (6) to get the score of the scheme: x = (0.9585, 0.9699, 0.9643, 0.9541, 0.9368).

(6) Sorting results of schemes : .x2 > x3 > x1 > x4 > x5. According to the above results, it is the best scheme: x2.

## 6.2 Comparative analysis

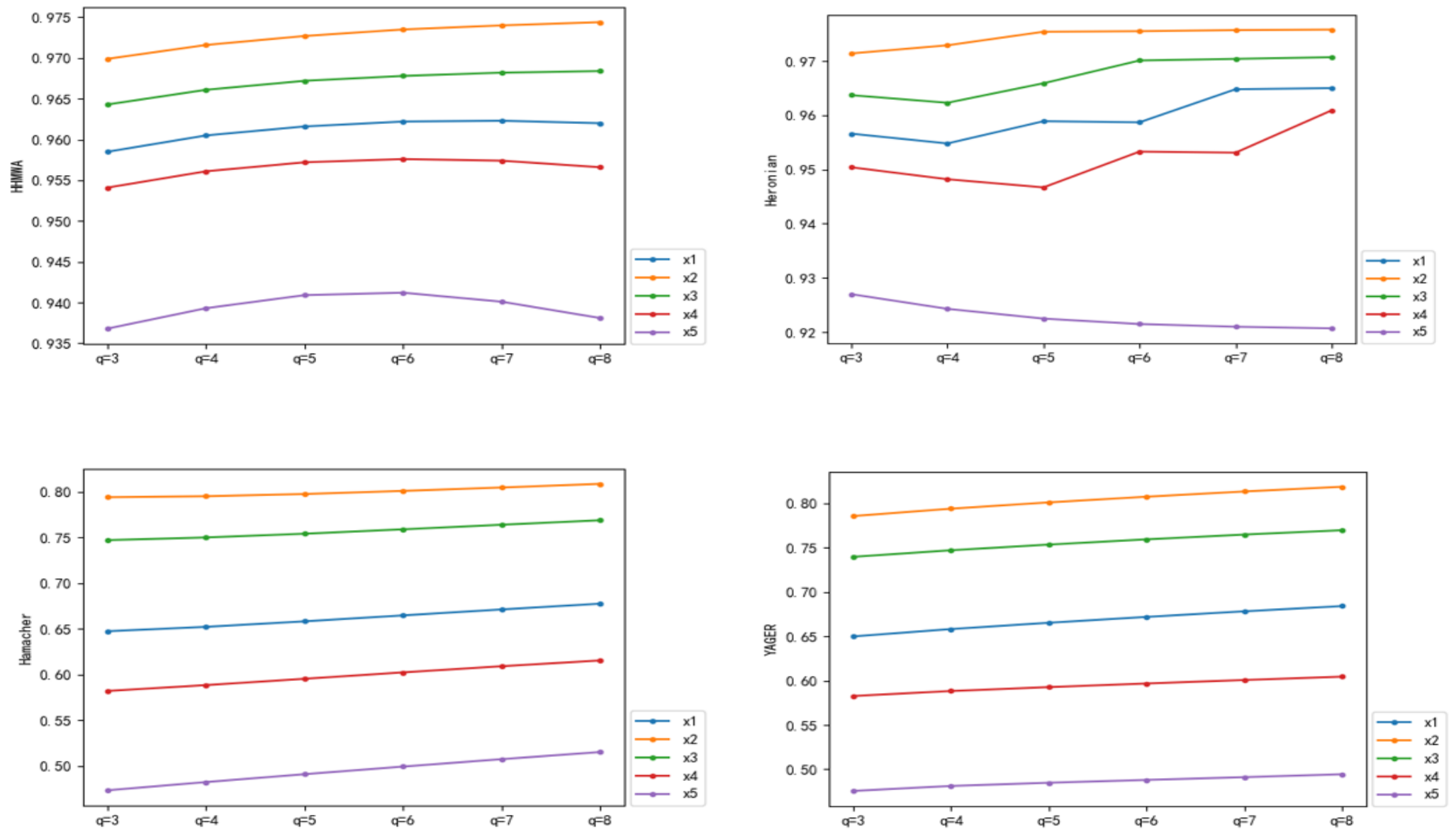

Figure 1 Comparative analysis results of operators

It can be seen from fig. 1 that the results of the proposed IVq-ROFHHMWA operator are in the same order as those of Haronian, Hamacher and Yager operators. And with the change of Q value, the ranking results of schemes are unchanged. Therefore, the validity and feasibility of the operator proposed in this paper are further verified.

# 7. Conclusion

In this paper, according to the basic rules of Hamacher operation, the algorithm based on interval valued q-rung Orthopair Fuzzy Hamacher operation is proposed, and then it is fused with Haronian operator, and the interval-valued generalized orthogonal HHM arithmetic average operator (IVq-ROFHHMWA) and interval valued q-rung Orthopair Fuzzy HHM geometric average operator (IVq-ROFHHMGA) are proposed. And their idempotent, boundary and permutation are studied. Then, based on IVq-ROFHHMWA operator, a multi-objective and multi-attribute group decision-making method is proposed. This method determines the value of Q according to the known data, and after data standardization transformation, the q-rung Orthopair Fuzzy numbers of each scheme after aggregation are obtained through two aggregations. Then, the aggregated Orthopair Fuzzy numbers are calculated and ranked by the score function and the exact function respectively, and the ranking results of each scheme and the optimal scheme can be obtained. Finally, the group decision-making method proposed in this paper is verified by a numerical example. The verification results show that the calculation results of the experimental case are consistent with the practical application and the operator operation results proposed by the existing experts, thus proving the feasibility and effectiveness of the operator and group decision-making method proposed in this paper.